\documentstyle[amstex,amscd,12pt]{article}
\makeatletter
\date{}
\newtheorem{theorem}{Theorem}[section]

\newtheorem{proposition}[theorem]{Proposition}

\newcommand{\edim}{{\rm e-dim}}
\newcommand{\z}{{\Bbb Z}}
\newcommand{\q}{{\Bbb Q}}

\newcommand{\invlim}{\raisebox{-1ex}{$\stackrel{\hbox{lim}}{\leftarrow}$}}

\newcommand{\lo}{\longrightarrow}
\newcommand{\sm}{\setminus}
\newcommand{\tor}{{\rm Tor}}
\begin{document}

\title{Acyclic resolutions for arbitrary groups\\
}

\author{Michael  Levin}

\maketitle
\begin{abstract}
We prove that for every abelian group $G$ and every compactum $X$ with
$\dim_G X  \leq n \geq 2$  there is a  $G$-acyclic resolution  $r: Z \lo X$ from a compactum $Z$
with  $\dim_G Z \leq n$ and  $\dim Z \leq n+1$  onto $X$.
\bigskip
\\
{\bf Math. Subj. Class.:} 55M10, 54F45.
\end{abstract}
\begin{section}{Introduction}

 Spaces denoted by $X$ are  assumed to be separable metrizable.  A compactum is
a metrizable compact space.

 Let $G$ be an abelian group.  A space   $X$  has the cohomological dimension
 $\dim_G X  \leq n$   if
 ${\check H}^{n+1}(X,A;G)=0$ for every closed
 subset $A$ of $X$. The case  $G=\z$ is an important special  case
 of cohomological dimension.  It  was known long ago that
 $\dim X  =\dim_\z X$ if $X$ is finite dimensional.  Solving an outstanding problem
 in cohomological dimension theory Dranishnikov constructed in 1987 an infinite
 dimensional compactum of $\dim_\z =3$. A few years earlier a  deep relation between
 $\dim_\z$ and $\dim$ was established by the  Edwards cell-like resolution
 theorem \cite{ed1,w1} saying that a compactum of $\dim_\z \leq n$ can be obtained
 as the  image of  a cell-like map defined  on a compactum of  $\dim \leq n$.
   A  compactum    $X$ is cell-like if any map $f : X \lo K$ from $X$ to a CW-complex  $K$  is
   null homotopic.
  A  map is cell-like if its fibers are cell-like. The reduced $\rm {\check C}$ech
 cohomology groups of  a cell-like compactum   are trivial with respect to any group $G$.

Acyclic resolutions  originated in the Edwards cell-like resolution.
A compactum $X$ is  $G$-acyclic if ${\Tilde{\Check H}}{}^* (X;G)=0$ and
a map is $G$-acyclic if its fibers are $G$-acyclic.  Thus a cell-like map is $G$-acyclic
with respect to any abelian group $G$.   By the Vietoris-Begle  theorem a $G$-acyclic map
cannot raise the cohomological dimension $\dim_G$.
Dranishnikov  proved the following important

 \begin{theorem}  {\rm (\cite{d0.5, d1})}
\label{t0}
Let $X$ be  a compactum with $\dim_\q X \leq n$, $ n \geq 2$.  Then there are a compactum
 $Z $  with $\dim_\q Z \leq n$
 and
 $\dim Z \leq n+1$ and a $\q$-acyclic map $r : Z \lo X$  from $Z$ onto $X$.
\end{theorem}

 It has been widely
conjectured that Theorem \ref{t0} holds for any abelian group $G$. A substantial  progress
in solving this conjecture was made by Koyama and  Yokoi \cite{ky2} who proved it for
a large class of groups including $\q$ and very recently by Rubin and Schapiro \cite{rs1}
who settled  the case
$G=\z_{p^\infty}$.

The purpose of this note is to finally answer  this conjecture affirmatively
by proving
  \begin{theorem}
\label{t1}
 Let $G$ be an abelian group and
 let $X$ be  a compactum with $\dim_G X \leq n$, $ n \geq 2$.  Then there are a compactum
 $Z $  with
 $\dim_G Z \leq n$ and  $\dim Z \leq n+1$ and a $G$-acyclic map $r : Z \lo X$  from $Z$ onto $X$.
\end{theorem}

In general the dimension $n+1$ of $Z$ in Theorem \ref{t1} is best possible \cite{ky2}.
However, it is unknown  if the dimension  of $Z$ in Theorem \ref{t0}
  can be reduced to  $n$.  In this connection let us also mention the
following interesting result of Dranishnikov.
  \begin{theorem} {\rm (\cite{d0})}
\label{d0}
 Let $X$ be  a compactum with $\dim_{\z_p} X \leq n$.  Then there are a compactum
 $Z $  with
  $\dim Z \leq n$ and a $\z_p$-acyclic map $r : Z \lo X$  from $Z$ onto $X$.
\end{theorem}
Our proof  of Theorem \ref{t1} essentially uses  Dranishnikov's idea of constructing
a $\q$-acyclic resolution presented in \cite{d1}
and involves  some methods of \cite{l1}.  The proof is self-contained and does not
rely on previous results concerning acyclic resolutions.       The paper  \cite{d1}  is
an excellent source of basic information on cohomological dimension theory.

 \end{section}
 \begin{section}{Preliminaries}

 All  groups below are  abelian  and functions between groups are homomorphisms.
$\cal P$ stands for the set of  primes.  For a non-empty subset $\cal A$ of $ \cal P$ let
$S({\cal A})=\{ p_{1}^{n_1}p_{2}^{n_2}...p_{k}^{n_k} : p_i \in {\cal A}, n_i \geq 0\}$
 be the set of  positive integers with prime factors from $\cal A$ and for the empty set define
  $S(\emptyset)=\{ 1 \}$.
Let   $G$ be a group and  $g \in G$.  We say that  $g$  is $\cal A$-torsion  if there is $n\in S(\cal A)$
such that $ng=0$ and   $g$  is $\cal A$-divisible if for every $n \in S(\cal A)$ there is $h\in G$
such that $nh=g$.
$\tor_{\cal A} G$ is the subgroup of the $\cal A$-torsion elements
of $G$.   $G$ is  $\cal A$-torsion  if $G=\tor_{\cal A}   G$, $G$ is $\cal A$-torsion free
if $\tor_{\cal A} G= 0$
and $G$ is $\cal A$-divisible if every element of $G$ is $\cal A$-divisible.

\begin{proposition}
\label{a1}
 ${}$

(i)   If $G$ is $\cal A$-torsion then $G$ is $(\cal P \sm \cal A)$-divisible
and  $(\cal P \sm \cal A)$-torsion free.

(ii) A factor group  of an $\cal A$-divisible
group is $\cal A$-divisible  and    a factor group of an $\cal A$-torsion group is $\cal A$-torsion.

(iii) The direct sum of $\cal A$-divisible groups is $\cal A$-divisible  and
      the direct sum of  $\cal A$-torsion  groups is $\cal A$-torsion.

Let $f: G\lo H$ be a homomorphism of groups $G$ and $H$ and
let  $H$  be  $\cal B$-torsion.
   Then $G /\tor_{\cal B} G$ is

(iv)    $\cal A$-divisible if
$\ker f$ is $\cal A$-divisible   and
$\cal B \subset \cal A$;

(v)    $\cal A$-torsion if
$\ker f$ is $\cal A$-torsion and
$\cal B \cap \cal A=\emptyset$;

(vi)    $\cal A$-torsion and $\cal A$-divisible
if    $\ker f$ is $\cal A$-torsion and $\cal A$-divisible and
 $\cal B \cap \cal A=\emptyset$.
\end{proposition}
  {\bf Proof.} The proof of {\bf (i), (ii), (iii)} is obvious.

   Let  $\phi :G \lo G/\tor_{\cal B}G$ be the projection and $\phi(x)=y$.
  Then there is $n\in S(\cal B)$ such that  $nf(x)=f(nx)=0$ and hence $nx \in \ker f$.

  {\bf (iv)}
   Let $m \in S(\cal A)$. Since ${\cal B} \subset {\cal A}$,
   $nm \in S({\cal A})$.
   Then  there is $z\in \ker f$
   such that $nmz=nx$. Hence $n(mz-x)=0$ and therefore
     $\phi(mz-x)=0$.  Thus $m\phi(z)=\phi(x)=y$ and $G/\tor_{\cal B} G$ is $\cal A$-divisible.

{\bf (v)}
    By (i) $\ker f$ is $(\cal P \sm \cal A)$-divisible
    and therefore
    there is $z\in \ker f$
   such that $nz=nx$.  Then $n(z-x)=0$ and there is $m \in S(\cal A)$ such that
   $mz=0$.  Hence $\phi(z)=\phi(x)=y$ and $my=\phi(mz)=0$ and
   (v) follows.

 {\bf    (vi)}
     By  (v)     $G /\tor_{\cal B} G$  is $\cal A$-torsion.
   By (i) $\ker f$ is $(\cal P \sm \cal A)$-divisible   and
since $\ker f$ is $\cal A$-divisible, $\ker f$ is $\cal P$-divisible. Then by (iv)
$G /\tor_{\cal B} G$  is $\cal A$-divisible.
 \hfill $\Box$
 \\
 \\

 The notation $\edim X \leq Y$ is used to indicate the property that  every map
 $f: A \lo Y$  of a closed subset $A$ of $X$ into $Y$ extends over $X$.  It is known
 that $\dim_G X \leq n$ if and only if $\edim X \leq K(G,n)$ where
 $K(G,n)$ is the Eilenberg-Mac Lane complex  of type $(G,n)$.
   A  map between CW-complexes  is  combinatorial if  the preimage of
   every subcomplex
   of the range is a subcomplex of the domain.

  Let $M$ be  a simplicial complex and let  $M^{[n]}$        be
    the $n$-skeleton of $M$ (=the union of all simplexes of $M$ of $\dim \leq n$).
By
 a resolution $EW(M,n)$   of $M$    we mean a CW-complex $EW(M,n)$ and
 a combinatorial map
 $\omega : EW(M,n) \lo  M$ such that $\omega$ is 1-to-1 over $M^{[n]}$.
 The resolution is said to be suitable for a map $f: M^{[n]} \lo  Y$   if
 the map  $f \circ\omega|_{\omega^{-1}(M^{[n]})}$ extends
   to a map from $ EW(M,n)$ to $ Y$. The resolution is said to be
   suitable  for  a compactum $X$
if for   every simplex $\Delta$ of $M$,
 $\edim X  \leq \omega^{-1}(\Delta)$.
     Note that if $\omega: EW(M,n) \lo M$ is a resolution suitable
 for $X$   then  for every map  $\phi :  X \lo  M$  there is  a map   $\psi : X \lo EW(M,n)$
 such that  for every simplex $\Delta$ of $M$ ,
  $(\omega \circ \psi)(\phi^{-1}(\Delta)) \subset \Delta$.
 We will call $\psi$ a combinatorial lifting of $\phi$.

  Following \cite{l1}  we  will construct  a   resolution of
 an $(n+1)$-dimensional simplicial complex $M$ which is
 suitable for  $X$ with $\dim_G X \leq n$  and a map $f : M^{[n]} \lo K(G,n)$.  In the sequel we will refer
 to this resolution as  the standard resolution for  $f$.  Fix  a CW-structure on
 $K(G,n)$    and assume that $f$ is cellular.
 We will obtain
 a CW-complex $EW(M,n)$   from  $M^{[n]}$ by attaching
   the mapping cylinder of  $f|_{\partial \Delta}$
 to      $\partial \Delta$  for every  $(n+1)$-simplex $\Delta$ of $M$.
 Let $\omega : EW(M,n) \lo M$ be the projection sending each mapping cylinder
 to  the corresponding $(n+1)$-simplex $\Delta$  such that $\omega$
 is the identity map on $\partial \Delta$,   the $K(G,n)$-part of the cylinder
 is sent to the barycenter  of     $\Delta$ and   $\omega$ is 1-to-1 on the rest of the cylinder.
 Clearly $f|_{\partial \Delta}$ extends over its mapping cylinder and therefore
   $f \circ\omega|_{\omega^{-1}(M^{[n]})}$  extends over $EW(M,n)$.
  For each simplex $\Delta$ of $M$,  $\omega^{-1}(\Delta)$ is either contractible or
  homotopy equivalent to $K(G,n)$.
    Define
   a CW-structure on      $EW(M,n)$   turning $\omega$ into a combinatorial map.
  Thus  we get  that
   the standard resolution is indeed a resolution suitable for both $X$ and $f$.
   Note that from the construction of the standard resolution $\omega : EW(M,n) \lo M$
   it follows that
    for every subcomplex $T$ of $M$, $\omega^{-1}(T)$ is the standard resolution of $T$
    for $f|_{T^{[n]}}$
   and    $\omega^{-1}(T)$
     is $(n-1)$-connected  if  $T$ is $(n-1)$-connected.

  \begin{proposition}
\label{h1}
Let $M$ be an $(n+1)$-dimensional finite simplicial complex and let
$\omega :EW(M,n) \lo M$ be
the  standard resolution for    $f: M^{[n]} \lo K(G,n)$, $n\geq 2$.  Then  for
$\omega_* : H_n(EW(M,n)) \lo  H_n(M) $,  $\ker \omega_*$  is a factor group of the
direct sum $\oplus G$ of finitely many $G$.
\end{proposition}
{\bf Proof.}   Inside each $(n+1)$-simplex of $M$ cut a small closed   ball around
the barycenter and not touching the boundary and split $M$ into two subspaces $M=M_1 \cup M_2$ where
$M_1$= the closure of
the complement  to the union of the balls and $M_2$=the union of   the balls.
Then $\omega$ is 1-to-1 over $M_1$,  $H_{n-1}(M_1\cap M_2) =0$,
$H_n(M_2)=0$  and the preimage under $\omega$ of each ball
is homotopy   equivalent to $K(G,n)$ and hence
$H_n(\omega^{-1}(M_2))$ is the direct sum $\oplus G$ of finitely many $G$.
  Consider  the Mayer-Vietoris
sequences for the pairs $(M_1,M_2)$ and $(\omega^{-1}(M_1), \omega^{-1}(M_2))$,
in which we  identify  $M_1$ and  $M_1\cap M_2$ with $\omega^{-1} M_1$  and
$\omega^{-1}(M_1\cap M_2)$  respectively.

 From the
       Mayer-Vietoris sequences
 it follows that
  $j_* (H_n(\omega^{-1}(M_1) \oplus H_n(\omega^{-1}(M_2)) )=H_n(\omega^{-1}(M_1 \cup M_2))$
  and $j_* (0\oplus  H_n(\omega^{-1}(M_2)) ) \subset \ker \omega_*$.
  Let us show that    $j_* (0\oplus  H_n(\omega^{-1}(M_2)) ) \supset \ker \omega_*$.
Let   $j_*(a\oplus b)\in \ker \omega_*$. Then in  the Mayer-Vietoris sequence
 for the pair $(M_1, M_2)$,  $j_*(a\oplus 0)=0$ and  therefore
 there is $c \in H_n(M_1\cap M_2))$ such that $i_*( c)=a\oplus 0$.
 Then in the Mayer-Vietoris sequence
 for the pair $(\omega^{-1}(M_1),\omega^{-1} ( M_2))$, $i_*(c)=a\oplus d $ and $j_* (a\oplus d)=0$.
 Thus $j_*(a\oplus b) = j_* (0 \oplus (b-d))$ and
 therefore  $ j_*(0\oplus  H_n(\omega^{-1}(M_2)))  = \ker \omega_*$.    Recall that
  $H_n(\omega^{-1}(M_2))=\oplus G$ and the proposition follows.
  \hfill $\Box$

   \begin{proposition}
\label{h2}
Let  $M=M_1 \cup M_2$  be a CW-complex with subcomplexes $M_1$ and $M_2$
such that  $M_1, M_2$ and $M_1 \cap M_2$  are
$(n-1)$-connected, $n\geq 2$ . Then $M$ is $(n-1)$-connected and

(i) $H_n(M)$ is $\cal A$-divisible if $H_n(M_1)$ and $H_n(M_2)$ are $\cal A$-divisible;

(ii)   $H_n(M)$ is $\cal A$-torsion  if $H_n(M_1)$ and $H_n(M_2)$ are $\cal A$-torsion.
\end{proposition}
 {\bf Proof.}  The connectedness of  $M$ follows from   van Kampen and  Hurewicz's theorems and
 the Mayer-Vietoris sequence.  (i) and (ii)
 follow from the Mayer-Vietoris sequence and (ii) and (iii) of
  Proposition \ref{a1}.
   \hfill $\Box$
   ${}$
  \\

  Let $X$ be a compactum  and let $\sigma (G)$ be  the Bockstein basis of
  a group  $G$.  By  Bockstein's  theory
   $\dim_G X \leq n$  if and only if  $\dim_E X \leq n$ for every $E \in \sigma (G)$.
    Denote :

   ${\cal T }(G)= \{ p \in{ \cal P} :  \z_p \in \sigma (G) \}$ ;

 ${\cal T }_\infty(G)= \{ p \in{ \cal P} :  \z_{p^\infty} \in \sigma (G) \}$ ;

  ${\cal D }(G)=  \cal P$ if $\q \in \sigma (G)$ and
  ${\cal D }(G)={\cal P} \sm  \{ p \in{ \cal P} :  \z_{(p)} \in \sigma (G) \}$       otherwise  ;

    ${\cal F}(G)={\cal D}(G) \sm ({\cal T}(G) \cup {\cal T}_\infty (G))$.

   Note that  ${\cal T }(G)$,    ${\cal T }_\infty(G)$
   and        ${\cal F}(G)$     are    disjoint  and $G$ is   ${\cal F}(G)$-torsion free.

 \begin{proposition}
\label{h3}
 Let $X$ be a compactum  and let $G$ be a group such that
 $G/{\tor G} \neq  0$ and $\dim_G X \leq n$. Then
   $\dim_E X \leq n$ for  every  group  $E$
   such that $E$ is ${\cal D}(G)$-divisible and ${\cal F}(G)$-torsion free.

\end{proposition}
 {\bf Proof.}   The proof  is based on Bockstein's theorem and inequalities.

 If $\z_p \in \sigma (E)$  then $\tor_p  E$ is not divisible by $p$ and hence  $E$ is not divisible by $p$.  Thus
 $p \in {\cal P} \sm{\cal D}(G)$ and therefore
  $\z_{(p)} \in  \sigma(G)$ and
 $\dim_{\z_p} X \leq \dim_{\z_{(p)} }X \leq n$.

 If    $\z_{p^\infty}  \in \sigma (E)$ then  $p$ is not in ${\cal F}(G)$.
 Then either   $p \in {\cal P} \sm{\cal D}(G)$
 and     $\dim_{\z_{p^\infty}} X \leq \dim_{\z_{(p)} }X \leq n$  or
  $p \in {\cal D}(G)\sm {\cal F}(G)$ and then either $p \in  {\cal T}(G)$
  and  $\dim_{\z_{p^\infty}} X \leq \dim_{\z_{p} }X \leq n$  or
  $p \in  {\cal T_\infty}(G)$ and  $\dim_{\z_{p^\infty}} X \leq n$.

  If $\z_{(p)} \in \sigma (E)$   then
  $E/{\tor E}$ is not divisible by $p$ and hence  $E$ is not divisible by $p$.  Thus
  $p \in {\cal P} \sm{\cal D}(G)$  and therefore
  $\dim_{\z_{(p)}} X  \leq n$.

  If  $\q \in  \sigma (E)$ then  consider the following cases:

   (i) ${\cal D}(G)=\cal P$.  Then since $G/\tor G \neq 0$, $\q \in \sigma(G)$
   and therefore $\dim_\q  X \leq n$ (this is the only place where we use
 that $G/\tor G \neq  0$);

  (ii)  there is      $p \in {\cal P} \sm{\cal D}(G)$. Then
   $\dim_\q X \leq \dim_{\z_{(p)}} X \leq n$.
   \hfill $\Box$

\end{section}
\begin{section}{Proof of Theorem \ref{t1}}
     Represent $X$ as   the inverse limit $X  = \invlim (K_i,h_i)$ of finite simplicial complexes $K_i$
 with combinatorial bonding maps   $h_{i+1}  :  K_{i+1} \lo K_i$ onto  and the projections
 $p_i : X \lo K_i$ such that for every simplex $\Delta$ of $K_i$,  diam$(p_i^{-1}(\Delta)) \leq 1/i$.
 Following A. Dranishnikov \cite{d1} we construct  by induction  finite CW-complexes
  $L_i$   and maps      $g_{i+1}: L_{i+1} \lo L_i$,
 $\alpha_i : L_i \lo K_i$   such that  \\

 (a)  $L_i$ is $(n+1)$-dimensional  and
  obtained from $K^{[n+1]}_i$ by replacing some $(n+1)$-simplexes by $(n+1)$-cells
 attached  to the   boundary of the  replaced simplexes by a map of degree $\in  S({\cal F}(G))$.
 Then $ \alpha_i$ is a projection of $L_i$ taking
 the new cells to the  original ones such that $\alpha_i$ is 1-to-1 over
 $K^{[n]}_i$.  We  define
 a simplicial structure on $L_i$ for which $\alpha_i$ is a  combinatorial map and refer
 to this simplicial structure while constructing resolutions of $L_i$.
 Note
 that for  ${\cal F}(G)=\emptyset$ we don't replace simplexes of $K^{[n+1]}_i$ at all;

 (b)  the maps $h_i$, $g_i$ and  $\alpha_i$   combinatorially commute.  By this we mean that
 for every simplex  $\Delta $ of $ K_{i}$,
   $(\alpha_i  \circ g_{i+1})((h_{i+1} \circ \alpha_{i+1})^{-1}(\Delta)) \subset \Delta$.
 \\

  We will construct $L_i$ in such a way that   $Z=\invlim(L_i,g_i)$ will be of $\dim_G \leq n$
  and $Z$  will admit a $G$-acyclic map  onto $X$.

 Set $L_1=K_1^{[n+1]}$ with $\alpha_1 : L_1 \lo K_1$ the embedding and assume that the construction
 is completed for $i$.   Let $E \in \sigma(G)$  and let $f  :  L_i^{[n]} \lo K(E,n)$ be a cellular map.
Let  $\omega_L: EW(L_i, n) \lo L_i^{}$
be the standard resolution  of $L_i$ for $f$.   We are going to construct from  $EW(L_i, n)$
a resolution of $K_i$ suitable for $X$.
On the first step of the construction we will
obtain from    $EW(L_i, n)$
   a resolution
 $\omega_{n+1} : EW(K_i^{[n+1]}, n) \lo K_i^{[n+1]}$ such that
 $EW(L_i, n)$ is a subcomplex of  $ EW(K_i^{[n+1]}, n)$  and $\omega_{n+1}$
 extends $\alpha_i \circ \omega_L$.   On the second
 step we will construct resolutions
 $\omega_j : EW(K_i^{[j]}, n) \lo K_i^{[j]}$, $n+2\leq j\leq \dim K_i$
  such that  $EW(K_i^{[j]}, n)$ is a subcomplex  of $EW(K_i^{[j+1]}, n)$
  and $\omega_{j+1}$ extends $\omega_j$. The construction  is carried out as follows.  \\

 Step 1.
For every simplex $\Delta$ of $K_i$ of $\dim =n+1$   consider
separately the subcomplex $(\alpha_i\circ\omega_L)^{-1}(\Delta)$
of $EW(L_i,n)$.
 Note that  from (a)  and the properties   of the standard resolution it follows
that the preimage under
   $\alpha_i\circ\omega_L $  of an $(n-1)$-connected subcomplex of $K_i$ is
   $(n-1)$-connected.
 Then  $(\alpha_i\circ\omega_L)^{-1}(\Delta)$
is $(n-1)$-connected. Enlarge  $(\alpha_i\circ\omega_L)^{-1}(\Delta)$
 by  attaching  cells  of $\dim = n+1$  in order to    kill
 $\tor_{{\cal F}(G)}H_n((\alpha_i\circ \omega_L)^{-1}(\Delta))$  and
  attaching cells of $\dim > n+1$
   in order to kill all homotopy groups of the enlarged subcomplex  in $\dim > n$.
 Define    $EW(K_i^{[n+1]}, n)$  as  $EW(L_i,n)$ with all the  cells attached for all
 $(n+1)$-dimensional simplexes $\Delta$ of $K_i$   and let
     a map $\omega_{n+1}:      EW(K_i^{[n+1]}, n)\lo K_i^{[n+1]}$        extend
     $\alpha_i \circ \omega_L$
     by sending the interior points of the attached cells to the interior of  the corresponding
    $\Delta$.  \\

 Step 2.  Assume that a  resolution   $\omega_j : EW(K_i^{[j]}, n) \lo K_i^{[j]}$,
 $n+1 \leq j < \dim K_i$  is constructed  such that the $n$-skeleton of   $ EW(K_i^{[j]}, n)$
 coincides   with the $n$-skeleton of  $EW(L_i, n)$.
 For every simplex $\Delta$ of $K_i$ of $\dim=j+1$
 consider separately  the subcomplex    ${\omega_j}^{-1}(\partial \Delta)$ of $EW(K_i^{[j]}, n)$.
 Then the $n$-skeleton of   ${\omega_j}^{-1}(\partial \Delta)$ coincides
 with the $n$-skeleton of     $(\alpha_i\circ\omega_L)^{-1}(\partial \Delta)$   and
 therefore      ${\omega_j}^{-1}(\partial \Delta)$ is $(n-1)$-connected.
 Enlarge ${\omega_j}^{-1}(\partial \Delta)$ by
  attaching  cells  of $\dim = n+1$ to
   in order to    kill
 $\tor_{{\cal F}(G)}H_n(\omega_j^{-1}(\partial \Delta))$  and
   attaching  cells of $\dim > n+1$        in order
  to kill all  homotopy groups of the enlarged  subcomplex  in $\dim > n$.
 Define    $EW(K_i^{[j+1]}, n)$  as  $EW(K_i^{[j]}, n) $ with all the cells
 attached
 for all $(j+1)$-simplexes of $K_i$  and let a map
     $\omega_{j+1}:      EW(K_i^{[j+1]}, n)\lo K_i^{[j+1]}$           extend $\omega_j$
     by sending the interior points of the attached cells to the interior of       the corresponding
    $\Delta$.   \\

    Finally denote  $EW(K_i, n)=EW(K_i^{[m]}, n)$ and $\omega=\omega_m : EW(K_i, n) \lo K_i$
    where $m=\dim K_i$.
From the construction
  it follows that the $n$-skeleton  of       $EW(K_i, n)$  is contained
  in    $EW(L_i,n)$     and
 for every
     simplex $\Delta$ of $K_i$,   $\omega^{-1}(\Delta )$ is contractible
      if $\dim \Delta \leq n$  and
       $\omega^{-1}(\Delta )$ is
    homotopy equivalent to $ K(H_n(  \omega^{-1}(\Delta)), n)$  if $\dim \Delta \geq n+1$.

    Let us show that $EW(K_i, n)$ is suitable for $X$.
 In order to verify that $\dim_{H_n(  \omega^{-1}(\Delta))}X \leq n$ for
 every simplex $\Delta $ of $K_i$   of  $\dim  \geq n+1$ we  first
    consider  Step 1 of
    the construction.  Let  $\Delta$ be  an $(n+1)$-dimensional simplex of $K_i$.
    By ${\omega_L}|_{...}$ we will denote  the map
    ${\omega_L}|_{( {\alpha_i}\circ {\omega_L})^{-1}({\Delta})} :
    ( {\alpha_i}\circ {\omega_L})^{-1}({\Delta}) \lo {\alpha_i}^{-1}(\Delta)$
    with the range restricted to ${\alpha_i}^{-1}(\Delta)$.
   Note that  by (a), $H_n(\alpha_i^{-1}(\Delta))$
    is ${\cal F}(G)$-torsion  and
    $H_n({\omega^{-1} }(\Delta))
    =H_n( (\alpha_i\circ \omega_L)^{-1}(\Delta))/\tor_{{\cal F}(G)}H_n( (\alpha_i
    \circ \omega_L)^{-1}(\Delta))$.
    Let   $( {\omega_L}|_{...})_*    :
 H_n(    ( {\alpha_i}\circ {\omega_L})^{-1}({\Delta}) )\lo H_n( {\alpha_i}^{-1}(\Delta))$.
 Consider the following cases.    \\

    Case 1-1.   $E=\z_p$.  By Proposition \ref{h1}
    $\ker( {\omega_L}|_{...})_*$ is
     $p$-torsion.
     Then  since $p$ is not in ${\cal F}(G)$, by  Proposition  \ref{a1}, (v),
       $H_n(\omega^{-1}(\Delta))$ is $p$-torsion and
       by Bockstein's theorem  $\dim_{H_n(  \omega^{-1}(\Delta))}X \leq  \dim_{\z_p} X  \leq n$.

    Case 1-2.  $E=\z_{p^\infty}$.  By Proposition \ref{h1}
    $\ker ({\omega_L}|_{...})_*$ is
     $p$-torsion  and $p$-divisible. Then  since $p$ is not in ${\cal F}(G)$, by  Proposition  \ref{a1}, (vi),
       $H_n(\omega^{-1}(\Delta))$ is $p$-torsion
       and $p$-divisible and  by  Bockstein's theorem
            $\dim_{H_n(  \omega^{-1}(\Delta))}X \leq \dim_{\z_{p^\infty}} X \leq n$.

    Case 1-3. $E=\z_{(p)}$ or $E=\q$.   By Proposition \ref{h1}
    $\ker ({\omega_L} |_{...})_*$
    is ${\cal D}(G)$-divisible.    Then  since ${\cal F}(G)\subset {\cal D}(G)$, by  Proposition  \ref{a1}, (iv),
       $H_n(\omega^{-1}(\Delta))$  is  ${\cal D}(G)$-divisible   and since
       $H_n(\omega^{-1}(\Delta))$
    is ${\cal F}(G)$-torsion free,
         by Proposition \ref{h3},
      $\dim_{H_n(  \omega^{-1}(\Delta))}X \leq n$.   \\

      Now let us pass to Step 2 of the construction.
      We will show that the properties of the homology
      groups  established  above will be preserved for
       simplexes of  higher dimensions. Let  $\Delta$ be
      a $(j+1)$-dimensional simplex of $K_i$,
      $j \geq n+1$ and recall that
    $H_n({\omega}^{-1}(\Delta))=
    H_n(\omega^{-1}(\partial \Delta))/\tor_{{\cal F}(G)}H_n(\omega^{-1}(\partial\Delta))$.
    Note that from the construction it follows that   the preimage under $\omega$ of
    an $(n-1)$-connected subcomplex of $K_i$  is $(n-1)$-connected.
    Also note that the intersection of a $j$-dimensional simplex
    of     $\Delta$      with the union of any collection of      $j$-dimensional simplexes
   of     $\Delta$         is     $(n-1)$-connected.
 These facts allow us to apply     below
     Proposition \ref{h2}  for  assembling $\omega^{-1}(\partial\Delta)$
     from  $\omega^{-1}(\Delta')$ for $j$-dimensional simplexes $\Delta'$ of $\Delta$
     to  show that   $\omega^{-1}(\partial\Delta)$
     has properties corresponding to  propperties of     $\omega^{-1}(\Delta')$.
    Once again we  consider separately the following cases.   \\

       Case 2-1.  $E=\z_{p}$.  If       for every $j$-dimensional
       simplex  ${\Delta}'$ of     $\Delta$,
          $H_n (\omega^{-1}({\Delta}'))$ is $p$-torsion then by Proposition \ref{h2},
         $    H_n( \omega^{-1}({\partial \Delta}))$     is      $p$-torsion and hence
   $ H_n({\omega}^{-1}(\Delta))$ is   $p$-torsion.
             Therefore   $\dim_{H_n(  \omega^{-1}(\Delta))}X \leq \dim_{\z_p} X \leq n$.

        Case  2-2.    $E=\z_{p^\infty}$.  If       for every $j$-dimensional simplex  ${\Delta}'$ of     $\Delta$,
          $H_n( \omega^{-1}({\Delta}'))$     is  $p$-torsion  and $p$-divisible then
           by Proposition \ref{h2},
        $    H_n( \omega^{-1}({\partial \Delta}))$   is         $p$-torsion  and $p$-divisible  and hence
         $ H_n({\omega}^{-1}(\Delta))$ is $p$-torsion  and $p$-divisible.
        Therefore
                 $\dim_{H_n(  \omega^{-1}(\Delta))}X \leq \dim_{\z_{p^{\infty}}}X \leq n$.

      Case 2-3. $E=\z_{(p)}$ or $E=\q$.               If       for every $j$-dimensional simplex  ${\Delta}'$ of     $\Delta$,
          $H_n( \omega^{-1}({\Delta}'))$   is   ${\cal D}(G)$-divisible then
            by Proposition \ref{h2},     $    H_n( \omega^{-1}({\partial \Delta}))$
            is      ${\cal D}(G)$-divisible. Then $ H_n({\omega}^{-1}(\Delta))$ is
            ${\cal D}(G)$-divisible and ${\cal F}(G)$-torsion free
              and    by Proposition \ref{h3},
      $\dim_{H_n(  \omega^{-1}(\Delta))}X \leq n$.  \\

         Thus we have shown that $EW(K_i, n)$ is suitable for $X$.  Now replacing $K_{i+1}$ by
         a $K_l$ with a sufficiently large $l$ we may assume that there is a combinatorial lifting
         of $h_{i+1}$ to $h'_{i+1} : K_{i+1} \lo  EW(K_i,n)$.   Replace $h'_{i+1}$ by its  cellular approximation
           preserving
  the property of   $h'_{i+1}$ of being a      combinatorial lifting of $h_{i+1}$.

   Consider the $(n+1)$-skeleton
 of $K_{i+1}$  and let $\Delta _{i+1}$ be an $(n+1)$-dimensional simplex
 in    $K_{i+1}$.   Let $\Delta_i $  be the smallest simplex  in $K_i$
 containing $h_{i+1}(\Delta_{i+1})$. Then $h'_{i+1}(\Delta_{i+1}) \subset\omega^{-1} (\Delta_i)$.
 Let ${\tau} : (\alpha_i\circ \omega_L)^{-1}(\Delta_i) \lo \omega^{-1} (\Delta_i)$
 be the inclusion. Note that  from the construction
it follows that for
 ${\tau}_* : H_n((\alpha_i\circ \omega_L)^{-1}(\Delta_i)) \lo H_n (\omega^{-1} (\Delta_i))$,
 $\ker \tau_*$ is ${\cal F}(G)$-torsion.
 Recall that the $n$-skeleton of $\omega^{-1} (\Delta_i)$ is contained in
$(\alpha_i\circ \omega_L)^{-1}(\Delta_i)$ and
 consider $h'_{i+1}|_{\partial \Delta_{i+1}} $ as a map
to  $(\alpha_i\circ \omega_L)^{-1}(\Delta_i)$. Let $a$ be the generator  of $H_n(\partial
\Delta_{i+1})$.  Since $h'_{i+1}|_{\partial \Delta_{i+1}} $ extends over $\Delta_{i+1}$ as a map
to $\omega^{-1} (\Delta_i)$ we have that $\tau_*((h'_{i+1}|_{\partial \Delta_{i+1}})_*(a))=0$.
Hence $(h'_{i+1}|_{\partial \Delta_{i+1}})_*(a) \in \ker \tau_*$ and therefore there is $k \in S
({\cal F}(G) )$ such that $ k((h'_{i+1}|_{\partial \Delta_{i+1}})_*(a))=0$.
  Replace $\Delta_{i+1}$ by a cell
     $C$  attached to the boundary of  $\Delta_{i+1}$ by a map of  degree   $k$
 if $(h'_{i+1}|_{\partial \Delta_{i+1}})_*(a)\neq 0$
     and set $C=\Delta_{i+1}$ if
$(h'_{i+1}|_{\partial \Delta_{i+1}})_*(a)=0$.
  Then $h'_{i+1}|_{\partial \Delta_{i+1}}$ can be
   extended over $C$  as a map to $(\alpha_i\circ \omega_L)^{-1}(\Delta_i)$
   and we will denote this extension by
 $g'_{i+1}|_{ C} : C\lo  (\alpha_i\circ \omega_L)^{-1}(\Delta_i)$.
 Thus replacing  if needed $(n+1)$-simplexes of $K_{i+1}^{[n+1]}$ we construct
 from $K_{i+1}^{[n+1]}$  a CW-complex $L_{i+1}$ and a map
 $g'_{i+1} :  L_{i+1} \lo  EW(L_i, n)$ which extends $h'_{i+1}$ restricted
 to  the $n$-skeleton  of $K_{i+1}$.   Now define
 $g_{i+1}=\omega_L \circ g'_{i+1} : L_{i+1} \lo L_i$ and
 finally define a simplicial structure on $L_{i+1}$ for which   $\alpha_{i+1}$
 is a combinatorial map.  It is easy to check that the properties (a) and (b)
 are satisfied.
  Since the triangulation of $L_{i+1}$ can be replaced
 by any of its barycentric subdivisions we  may also  assume that  \\

 (c)
 diam$g_{i+1}^j(\Delta) \leq 1/i$ for every simplex
 $\Delta$ in  $L_{i+1}$  and $j \leq i$   \\
     where
 $g^j_i=g_{j+1} \circ g_{j+2} \circ ...\circ g_i: L_i \lo L_j$.   \\

           Denote $Z=\invlim(L_i,g_i)$ and let $r_i : Z \lo L_i$ be the projections.
            For constructing
           $L_{i+1}$  we used an arbitrary map $f : L_i^{[n]} \lo K(E,n), E\in \sigma(G)$.
            Let us show that  choosing $E \in \sigma(G)$ and $f$  in an appropriate way
           for each $i$ we can achieve
           that $\dim_E Z \leq n$  for every  $   E \in \sigma(G)$ and hence  $\dim_G Z \leq n$.

  Let $\psi : F \lo K(E,n)$ be a map of a closed subset
 $F$ of $L_j$.
  Then by (c) for a sufficiently large $i>j$ the map
 $\psi \circ g_i^j|_{(g_i^j)^{-1}(F)}$ extends over a subcomplex $N$ of $L_i$
to  a map $\phi : N \lo  K(E,n)$.  Extending
    $\phi$ over  $L_i^{[n]}$ we may assume that $   L_i^{[n]} \subset N$ and replacing
   $\phi$ by its cellular approximation we  assume that      $\phi$
   is cellular.
    Now define  the map  $f :L_i^{[n]} \lo K(E,n)$ that we use  for
    constructing $L_{i+1}$  as   $f =\phi |_{ L_i^{[n]}} $.
          Since $g_{i+1}$ factors
           through $  EW(L_i, n)$, the   map $f \circ g_{i+1}|_{g^{-1}_{i+1}(L^{[n]}_i)}:
 {g^{-1}_{i+1}(L^{[n]}_i)} \lo K(E, n)$ extends to a map $f' : L_{i+1}\lo K(E, n)$.
  Define $\psi ' :  L_{i+1}\lo K(E, n)$ by $\psi'(x)=(\phi \circ g_{i+1})(x)$
  if $x \in {g_{i+1}^{-1}}(N)$ and $\psi'(x)=f'(x)$ otherwise.
  Then $\psi'|_{(g^{j}_{i+1})^{-1}(F)} : {(g^{j}_{i+1})^{-1}(F)} \lo K(E, n)$
   is homotopic to
  $\psi \circ g_{i+1}^j |_{(g^{j}_{i+1})^{-1}(F)}:
 {(g^{j}_{i+1})^{-1}(F)} \lo K(E, n)$ and hence
 $\psi \circ g_{i+1}^j |_{(g^{j}_{i+1})^{-1}(F)}$ extends over $L_{i+1}$.
Now since we need to solve only countably many extension problems
for every $L_j$ with respect to $K(E,n)$ for every $E \in \sigma(G)$
we can choose
  for each $i$ a map $f: L_i^{[n]} \lo  K(E,n)$   in the way
  described above to achieve that  $\dim_E Z \leq n$ for every $E \in \sigma(G)$  and
           hence $\dim_G Z \leq n$.  \\

           The property  (b) implies that
           for every   $x \in X$ and  $z\in Z$,

 (d1) $ g_{i+1}( \alpha_{i+1}^{-1}(st( p_{i+1} (x) ) )) \subset
   \alpha_{i}^{-1}(st( p_i (x)) )$ and

     (d2)
   $h_{i+1}(st((\alpha_{i+1}\circ r_{i+1})(z)))  \subset        st((\alpha_{i}\circ r_{i})(z))$
 \\ where $st(a)=$the union of all the simplexes containing $a$.  \\

   Define a  map $r : Z \lo X$ by $r(z)=\cap \{ p_i^{-1}(    st((\alpha_{i}\circ r_{i})(z) ) ): i=1,2,... \}$.
  Then  (d1)  and (d2) imply
   that $r$ is indeed well-defined and continuous.

  The properties (d1) and (d2) also imply  that  for every $x \in X$
\\ $r^{-1}(x)=\invlim ( \alpha_i ^{-1}(st(p_{i} (x))), g_i |_{\alpha_i ^{-1}(st(p_{i} (x)))})$,
   where  the map  $ g_i |_{\alpha_i ^{-1}(st(p_{i} (x)))}$ is considered as a map
      to $\alpha_{i-1} ^{-1}(st(p_{i-1} (x)))$.

Since $r^{-1}(x) $ is not empty for every $x \in X$,
 $r$ is a map onto and let us show that $r^{-1}(x)$ is
$G$-acyclic.

    Since $  st(p_{i} (x))$ is contractible,   $T=\alpha_i ^{-1}(st(p_{i} (x))) $ is $(n-1)$-connected.
          From (a) and Proposition \ref{h2}
   it follows that   $H_n(   T )$ is ${\cal F}(G)$-torsion.
Then, since $G$ is ${\cal F}(G)$-torsion free,
by  the universal-coefficient theorem $H^n(T;G)={\rm Hom}(H_n(T),G)=0$.
    Thus     ${\Tilde {\Check H}}{}^k(r^{-1}(x); G) =0$ for
     $ k\leq n$ and
 since $\dim_G Z \leq n$,
    ${\Tilde  {\Check H}}{}^k(r^{-1}(x); G) =0$ for  $ k\geq n+1$.
      Hence
  $r$ is $G$-acyclic and  this completes the proof.
  \hfill $\Box$

\end{section}

Department of Mathematics\\
Ben Gurion University of the Negev\\
P.O.B. 653\\
Be'er Sheva 84105, ISRAEL  \\
e-mail: mlevine@math.bgu.ac.il\\\\
\end{document}